\documentclass[10pt]{amsart}
\usepackage{a4}
\usepackage{amssymb}
\usepackage{setspace}
\newcommand{\frow}{{}^\frown }

\newcommand{\rc}{\operatorname{rc}}
\newcommand{\lpr}{\operatorname{lpr}}
\newcommand{\lex}{\operatorname{lex}}
\newcommand{\dom}{\operatorname{dom}}
\newcommand{\Fr}{\operatorname{Fr}}

\newcommand{\free}{\operatorname{free}}
\newcommand{\id}{\operatorname{id}}
\newcommand{\cf}{\operatorname{cf}}
\newcommand{\NPT}{\operatorname{NPT}}
\newcommand{\NRT}{\operatorname{NRT}}
\newcommand{\card}[1]{\mid\!\!{#1}\!\!\mid}
\newtheorem{thm}{Theorem}[section]
\theoremstyle{definition}
\newtheorem{defn}[thm]{Definition}
\newtheorem{lemma}[thm]{Lemma}
\newtheorem{question}[thm]{Question}
\newtheorem{corollary}[thm]{Corollary}

\begin{document}

\title{The number of openly generated Boolean algebras}

\author{Stefan Geschke}
\address[Geschke]{II.~Mathematisches Institut\\Freie Universit\"at
Berlin\\Arnimallee 3\\14195 Berlin\\Germany}
\email{geschke@math.fu-berlin.de}
\author{Saharon Shelah}
\address[Shelah]{Institute of Mathematics, The Hebrew University of 
Je\-ru\-sa\-lem, 91904 Je\-ru\-sa\-lem, Israel and
Department of Mathematics, Rutgers University, New Brunswick, NJ 
08854, USA.}
\email{shelah@math.huji.ac.il}
\date{March 19, 2004.}
\subjclass[2000]{06E05}
\keywords{projective Boolean algebra, openly generated, almost free}
\thanks{This article is [GeSh:558] in the second author's list of 
publications.  
The research of both authors was at least partially supported
by the Edmund Landau Center for research
in Mathematical Analysis, supported by the Minerva Foundation
(Germany).}

\begin{abstract} 
This article is devoted to two different generalizations of projective Boolean
algebras:  
openly generated Boolean algebras and tightly $\sigma$-filtered
Boolean algebras.  

We show that for every uncountable regular cardinal $\kappa$ there are $2^\kappa$ pairwise non-isomorphic openly generated 
Boolean algebras of size $\kappa>\aleph_1$ provided there is an
almost free non-free abelian group of size $\kappa$.
The openly generated Boolean algebras constructed here
are almost free.  
  
Moreover, for every infinite regular cardinal $\kappa$ we construct 
$2^\kappa$ pairwise non-isomorphic Boolean algebras of size $\kappa$ that
are tightly 
$\sigma$-filtered and c.c.c.

These two results contrast nicely with Koppelberg's
theorem in \cite{koppproj} 
that for every uncountable regular cardinal $\kappa$ there are only $2^{<\kappa}$ 
isomorphism types of projective Boolean algebras of size $\kappa$. 
\end{abstract}

\maketitle

\section{introduction}
Projectivity is usually defined as a universal property.  
A Boolean algebra $A$ is projective if and only if for every Boolean
algebra
$B$ and every epimorphism $f:B\to A$ there is a homomorphism $g:A\to B$
such that $f\circ g=\id_A$.  
However, theorems of Haydon, Koppelberg, and \v{S}\v{c}epin provide an
internal characterization of projectivity for Boolean algebras (see
\cite{koppproj}).  

Using her characterization of projectivity,  Koppelberg \cite{koppproj} 
showed that for every uncountable cardinal $\kappa$ there are only
$2^{<\kappa}$ isomorphism types of projective Boolean algebras of size
$\kappa$.  
She also showed that for every singular cardinal $\mu$ there
are $2^\mu$ pairwise non-isomorphic Boolean algebras of size $\mu$. 

For Boolean algebras, there are two natural generalizations of
projectivity: open generatedness (or rc-filteredness) and tight
$\sigma$-filteredness. 

Openly generated Boolean algebras are studied to a great
extend in \cite{heinsha}.   We only note that openly generated Boolean
algebras seem to be
quite close to projective Boolean algebras.  Every openly generated
Boolean algebra of size $\leq\aleph_1$ is projective and it is a non-trivial
task to find an openly generated Boolean algebra which is not projective. 
Examples of openly generated, non-projective Boolean algebras
were provided by \v{S}\v{c}epin (see \cite{heinsha}).  

Tight $\sigma$-filteredness was introduced in \cite{kopptight} 
and studied systematically in \cite{ge2}.   
It turns out that several properties of projective Boolean algebras 
can be generalized to tightly $\sigma$-filtered Boolean algebras.  

However, in the present article we show that for each infinite regular
cardinal $\kappa$ there are $2^\kappa$ pairwise non-isomorphic tightly
$\sigma$-filtered Boolean algebras of size $\kappa$, contrasting
Koppelberg's result on the number of projective Boolean algebras.  
The construction is fairly elementary and even yields Boolean algebras
which are c.c.c.   
This result is contained in the second author's PhD thesis \cite{ge}.

With openly generated Boolean algebras the situation is more complicated. 
Fu\-chi\-no  
showed that if $\kappa$ is a regular cardinal that has a 
non-reflecting  stationary subset consisting of ordinals of cofinality $\aleph_1$, 
then there are $2^\kappa$ pairwise non-isomorphic openly generated Boolean 
algebras of size $\kappa$ (see \cite{heinsha}).
Fuchino's assumption, i.e.,
the existence of a non-reflecting stationary set  
of ordinals of cofinality $\aleph_1$, implies the existence of 
an almost free, non-free abelian group of size $\kappa$ 
\cite[Chapter VI, Lemma 2.2 and Theorem 2.3]{ekmek}.  
Here an abelian group is almost free if every subgroup of
strictly smaller cardinality is free.  

In the present article we show that for an infinite regular cardinal
$\kappa$ there are $2^\kappa$ pairwise non-isomorphic openly generated
Boolean
algebras if there exists an almost free, non-free abelian group of size
$\kappa$.

While every almost free abelian group of singular cardinality is in
fact free by Shelah's
compactness theorem (see \cite{ekmek}), there are almost free, non-free 
abelian groups of various cardinalities.  
The first example of such a group is due to Baer (see \cite{fuchs}) and
has cardinality $\aleph_1$.
Let $\mathcal C$ denote the class of cardinals $\kappa$ for
which there is an almost free, non-free abelian group of size $\kappa$.  
By Baer's result, $\aleph_1\in\mathcal C$.
  
In \cite{mashe} Magidor and Shelah showed that
$\mathcal C$ is closed
under taking successors and under the operation
$(\lambda,\kappa)\mapsto\lambda^{+(\kappa+1)}$. 
On the other hand, they proved, assuming the consistency of infinitely
many supercompact cardinals, that it is consistent that $\mathcal C$
contains no cardinal above the first cardinal fixed point, i.e., the first
$\kappa$ with $\kappa=\aleph_\kappa$.
However, if $V=L$, then every regular cardinal $\kappa$ which
is not weakly compact belongs to $\mathcal C$.  A proof of the latter fact
and much more on this topic can be found in the book
\cite{ekmek} by Eklof and Mekler.   
Note that there is a new, revised edition of this book \cite{ekmek2}.

\section{Basic definitions and preliminary Lemmas}
Open generatedness and tight $\sigma$-filteredness are both defined in
terms of nicely embedded subalgebras.  Openly generated Boolean algebras
have many relatively complete subalgebras and tightly $\sigma$-filtered
Boolean algebras have many $\sigma$-sub\-al\-ge\-bras.  

\begin{defn} 
Let $A$ and $B$ be Boolean algebras such that $A\leq B$.
Then for $b\in B$ the ideal $\{a\in A:a\leq b\}$ of $A$ is denoted as
$A\restriction b$.  $A$ is called a {\em relatively complete
subalgebra} ({\em rc-subalgebra})
of $B$ if  for each $b\in B$ the ideal $A\restriction b$ is
principal.  In this case we write $A\leq_{\rc}B$.  If $A\leq_{\rc}
B$, then
$\lpr^B_A$ denotes the mapping from $B$ to $A$ assigning to each $b\in B$
the generator of $A\restriction b$, the lower projection of $b$ in $A$.

$A$ is called a {\em $\sigma$-subalgebra} of $B$ if for every $x\in B$
the ideal $A\restriction x$ is countably generated.
In this case we write $A\leq_\sigma B$.
If $C$ is a Boolean algebra and $e:C\to B$ an embedding, then $e$ is an 
{\em rc-embedding} if $e[C]\leq_{\rc}B$ and an 
{\em $\sigma$-embedding} if $e[C]\leq_\sigma B$.
\end{defn}

Note that $A\leq_\sigma B$ if and only if $A\leq B$ and for every ideal
$I$ of $B$ which is countably generated, $I\cap A$ is countably generated
as well.

\begin{defn} 
Let $A$ be a Boolean algebra and $\delta$ an ordinal. 
A sequence $(A_\alpha)_{\alpha<\delta}$ of subalgebras of $A$ is a {\em
filtration} of $A$ if 
\begin{itemize} 
\item[a)] $\bigcup_{\alpha<\delta}A_\alpha=A$,
\item[b)] $A_\alpha\leq A_\beta$ for $\alpha<\beta<\delta$, and
\item[c)] the sequence $(A_\alpha)_{\alpha<\delta}$ is continuous, i.e.,
for every limit ordinal
$\beta<\delta$, $A_\beta=\bigcup_{\alpha<\beta}A_\alpha$. 
\end{itemize}

A filtration $(A_\alpha)_{\alpha<\delta}$ of $A$ is {\em tight}
if there is a sequence $(x_\alpha)_{\alpha<\delta}$ in $A$ such that 
for all $\beta<\delta$, $A_\beta=\langle\{x_\alpha:\alpha<\beta\}\rangle$.
Here for a set $X\subseteq A$, $\langle X\rangle$ denotes the subalgebra
of $A$ generated by $X$.   
If $B\leq A$ and $x\in A$, we write
$B(x)$ for $\langle B\cup\{x\}\rangle$.  

$A$ is tightly $\sigma$-filtered if it has a tight filtration
$(A_\alpha)_{\alpha<\delta}$ such that for all $\alpha<\delta$,
$A_\alpha\leq_\sigma A$.  $(A_\alpha)_{\alpha<\delta}$ is called a {\em
tight $\sigma$-filtration} of $A$.

$A$ is openly generated if the set of rc-subalgebras of $A$ includes a club
of $[A]^{\leq\aleph_0}$.  
\end{defn}

Koppelberg's characterization of projective Boolean algebras is obtained 
by replacing $\leq_\sigma$ by $\leq_{\rc}$ in the definition of tight
$\sigma$-filteredness. 
The term ``openly generated'' comes from the fact that rc-embeddings 
correspond to open mappings via Stone duality. 

Why are projective algebras openly generated?  
For a set $X$ let $\Fr(X)$ denote the free Boolean algebra over the set
$X$. If $X\subseteq Y$, we regard $\Fr(X)$ as a subalgebra of $\Fr(Y)$ in
the obvious way.  It is easy to see that $\Fr(X)\leq_{\rc}\Fr(Y)$ for 
$X\subseteq Y$.  It follows that free Boolean algebras are openly
generated.  

By abstract nonsense,  
a Boolean algebra $A$ is projective if and only if it is a retract of a
free Boolean algebra, i.e., if there are a free Boolean algebra $B$ and
homomorphisms $e:A\to B$ and $p:B\to A$ such that $p\circ e=\id_A$.  
It is more or less straight forward to see that open generatedness is
hereditary with respect to retracts.  It follows that projective Boolean
algebras are openly generated.  

We collect some facts on rc-embeddings, $\sigma$-embeddings, and open
generatedness.

\begin{lemma} \label{simple}
Let $A\leq B$ and $x\in B$. Then the following hold:
\begin{itemize}\item[a)] $A\leq_{\operatorname{rc}} A(x)$ if and only if
$A\restriction x$ and $A\restriction -x$ are principal.
\item[b)] $A\leq_\sigma A(x)$ if and only if $A\restriction x$ and
$A\restriction -x$ are countably generated.
\end{itemize}
\end{lemma}

\begin{proof} We show b) only.  The direction from the left to the right is
trivial.
For the other direction let
$E\subseteq A$ and $F\subseteq A$ be countable sets which generate 
$A\restriction x$ and $A\restriction -x$ respectively.  Without
loss of generality we may
assume that $E$
and $F$ are closed under finite joins.  Suppose $y\in A(x)$.  Then 
there are $v,w\in A$ s.t.~$y=(v+x)\cdot(w+(-x))$.  Let $z\in A$ such that
$z\leq
y$. Then
$z-v\leq x$ and $z-w\leq -x$.  Hence $z-v\leq a$ and $z-w\leq b$ for some
$a\in E$ and some $b\in F$. It follows that $z\leq (v+a)\cdot(w+b)$.
Clearly, $(v+a)\cdot(w+b)\leq y$ for every $a\in E$ and every $b\in F$.
Hence  $A\restriction y$ is generated by $\{(v+a)\cdot(w+b):a\in E\wedge
b\in F\}$.
\end{proof}

\begin{lemma}
\label{trans}
\begin{itemize}\item[a)] 
$A\leq_{\rc} B\leq_{\rc} C$ $\Rightarrow$ $A\leq_{\rc} C$,
\item[] $A\leq_\sigma B\leq_\sigma C$ $\Rightarrow$ $A\leq_\sigma C$.
\item[b)] $A\leq B$, $B=\bigcup_{\alpha\leq\lambda}B_\alpha$, $A\leq_{\rc} 
B_\alpha$ for every $\alpha<\lambda$ $\Rightarrow$ $A\leq_{\rc} B$,
\item[] $A\leq B$, $B=\bigcup_{\alpha\leq\lambda}B_\alpha$,
$A\leq_\sigma  B_\alpha$ for every $\alpha<\lambda$ $\Rightarrow$
$A\leq_\sigma B$.
\item[c)] $A\leq_{\rc}B$, $C\leq B$ and $\lpr_A^B[C]\subseteq C$ $\Rightarrow$ 
$A\cap C\leq_{\rc}C$. 
\end{itemize}
\end{lemma}

\begin{proof} Easy. \end{proof}

\begin{lemma}\cite[Proposition 2.2.4]{heinsha}\label{provingopen} 
If $\delta$ is an ordinal and $A$ is
the union of an increasing continuous chain 
$(A_\alpha)_{\alpha<\delta}$ of rc-subalgebras that are openly generated,
then $A$ itself is openly generated.  
\end{lemma}

\section{The Number of Tightly $\sigma$-Filtered Boolean Algebras}

We show that Koppelberg's result on the number of projective Boolean
algebras cannot be extended to tightly $\sigma$-filtered Boolean algebras.

\begin{thm}\label{tightmany} For every infinite cardinal $\kappa$ there
are $2^\lambda$
pairwise non-isomorphic tightly
$\sigma$-filtered c.c.c.~Boolean algebras of size $\kappa$. \qed
\end{thm}  

The proof of the theorem uses the following lemma, which says that
stationary 
sets consisting of ordinals of countable cofinality can be coded by tightly
$\sigma$-filtered Boolean algebras.

\begin{lemma}\label{tightcoding}
Let $\kappa$ be an uncountable regular cardinal and
let $S$ be a subset of $\kappa$ consisting of ordinals of cofinality
$\aleph_0$.
Then there are a Boolean algebra $A$ of size $\lambda$ and a tight 
$\sigma$-filtration 
$(A_\alpha)_{\alpha<\kappa}$ of $A$ such that the following hold:
\begin{itemize}\item[a)] $A_\alpha\not\leq_{\rc}A$ for all $\alpha\in S$
\item[b)] $A_\alpha\leq_{\rc}A$ for all $\alpha\in\kappa\setminus S$.
\end{itemize}
\end{lemma}

\begin{proof}
For every $\alpha\in S$ let $(\delta_n^\alpha)_{n\in\omega}$ be a strictly 
increasing sequence of ordinals with least upper bound $\alpha$ and
$S\cap\{\delta_n^\alpha:n\in\omega\}=\emptyset$.  We will construct
$(A_\alpha)_{\alpha<\kappa}$ together with a sequence 
$(x_\alpha)_{\alpha<\kappa}$ such that
\begin{itemize}\item[(i)] $A_0=2$,
\item[(ii)] $A_{\alpha+1}=A_\alpha(x_\alpha)$ for all $\alpha<\kappa$,
\item[(iii)] $x_\alpha$ is independent over $A_\alpha$ whenever 
$\alpha\not\in S$,
\item[(iv)] $A_\alpha\restriction x_\alpha$ is generated by
$\{x_{\delta_n^\alpha}:n\in\omega\}$ and 
$A_\alpha\restriction-x_\alpha=\{0\}$ whenever
$\alpha\in S$, 
\item[(v)] $A_\beta=\bigcup_{\alpha<\beta}A_\alpha$ holds for all limit 
ordinals $\beta<\kappa$.
\end{itemize}
Clearly, the construction can be done and is uniquely determined.
We have to show that a) and b) of the lemma hold for 
$(A_\alpha)_{\alpha<\kappa}$.

For a) let $\alpha\in S$. Then $A_\alpha\restriction x_\alpha$ 
is non-principal.  For suppose
$a\in A_\alpha$ is such that $a\leq x_\alpha$.  
Since $(\delta_n^\alpha)_{n\in\omega}$ is
cofinal in $\alpha$, there is $n\in\omega$ 
such that $a\in A_{\delta_n^\alpha}$.  Since
$\delta_n^\alpha\not\in S$, $x_{\delta_n^\alpha}$ is independent over
$A_{\delta_n^\alpha}$ by construction.
Hence $a+x_{\delta_n^\alpha}$ is strictly larger than $a$,
but still smaller than
$x_\alpha$.  So $a$ does not generate $A_\alpha\restriction x_\alpha$.

For b) let $\alpha\not\in S$.  By induction on $\gamma<\kappa$, we show
that
$A_\alpha\leq_{\rc}A_\gamma$ holds for every $\gamma\geq\alpha$.
$A_\alpha\leq_{\rc}A_\alpha$ holds trivially.  
Suppose $\gamma$ is a limit 
ordinal and $A_\alpha\geq_{\rc}A_\beta$ holds for all $\beta<\gamma$
such that $\alpha\leq\beta$.  
Then $A_\alpha\leq_{\rc}A_\gamma$ follows from Lemma \ref{trans}.
Now suppose $\gamma=\beta+1$ for some $\beta\geq\alpha$.
There are two cases:

\begin{itemize}
\item[I.] $\beta\not\in S$.  In this case $A_\beta\leq_{\rc}A_\gamma$ by
construction.  
By hypothesis, $A_\alpha\leq_{\rc}A_\beta$.  By Lemma \ref{trans},
this implies $A_\alpha\leq_{\rc}A_\gamma$.

\item[II.] $\beta\in S$.  This is the non-trivial case.   We {\bf claim}
that
$A_\delta\leq_{\rc}A_\delta(x_\beta)$  holds for every $\delta<\beta$.
This can be seen as follows:  By Lemma \ref{simple}, it is sufficient to
show that both
$A_\delta\restriction x_\beta$ and $A_\delta\restriction -x_\beta$ 
are principal.  But
$A_\delta\restriction -x_\beta\subseteq A_\beta\restriction -x_\beta=\{0\}$ 
by construction.
Let $a\in A_\delta$ be such that $a\leq x_\beta$.  Let
$m:=\{n\in\omega:x_{\delta_n^\beta}\in A_\delta\}$.  Clearly $m\in\omega$.  
Let
$T\in[\omega]^{<\aleph_o}$ be such that
$a\leq\sum\{x_{\delta_n^\beta}:n\in T\}$. 
Then $$a\leq\sum\{x_{\delta_n^\beta}:
{n\in T\cap m}\}+\sum\{x_{\delta_n^\beta}:{n\in
T\setminus m}\}.$$ 
Since $\sum\{x_{\delta_n^\beta}:{n\in T\setminus m}\}$ is independent over 
$A_\delta$ by construction,  
$$a\leq\sum\{x_{\delta_n^\beta}:{n\in T\cap
m}\}\leq\sum\{x_{\delta_n^\beta}:{n<m}\}\leq
x_\beta.$$   This shows that $A_\delta\restriction x_\beta$ is generated by 
$\sum\{x_{\delta_n^\beta}:{n<m}\}$ and the claim follows.

Now 
$A_\gamma=A_\beta(x_\beta)=\bigcup_{\alpha\leq\delta<\beta}A_\delta(x_\beta)$.
Hence, $A_\alpha\leq_{\rc}A_\gamma$ follows from the claim together with
Lemma \ref{trans}.
\end{itemize}
This shows b).
\end{proof}

In order to show that the Boolean algebra $A$ constructed in the lemma above 
is c.c.c., we use an argument which was used 
in an early version of \cite{fugesoshe} to prove
that, assuming the consistency of the existence of a supercompact
cardinal,
it is consistent with ZFC+GCH
that there is a complete c.c.c.~Boolean algebra without the so-called WFN, 
a property whose definition is obtained by replacing $\leq_{\rc}$ by
$\leq_\sigma$ and $[A]^{\leq\aleph_0}$ by $[A]^{\leq\aleph_1}$ in the
definition of open generatedness.

\begin{lemma}\label{ccc}
The Boolean algebra $A$ constructed in the proof of 
Lemma \ref{tightcoding} is c.c.c.
\end{lemma}

\begin{proof}
Assume $A$ is not c.c.c.  Let $C\subseteq A$ be an
uncountable
antichain.  Let $X:=\{x_\alpha:\alpha<\kappa\}$.  For $x\in X$ let
$x^0:=x$ and $x^1:=-x$.  
We may assume that each $a\in C$ is an elementary product of elements of
$X$, 
i.e., there is
$X_a\in[X]^{<\aleph_0}$ and $f_a:X_a\to 2$ such that 
$a=\prod_{x\in X_a}x^{f_a(x)}$.
After thinning out $C$ if necessary,
we may assume that $\{X_a:a\in C\}$ is a $\Delta$-system with root $R$, 
there is $f:R\to 2$ such that $f_a\restriction R=f$ for all $a\in C$, and all
$X_a$ are of the same size, say $n$.


{\bf Claim.} Let $Y\in [X]^{<\omega}$ and $g:Y\to 2$ be such that
$\prod_{x\in
Y}x^{g(x)}=0$. Then there are $\alpha\in S$ and $i\in\omega$ with
$x_\alpha,x_{\delta^\alpha_i}\in Y$ such that $g(x_\alpha)=1$ and
$g(x_{\delta^\alpha_i})=0$.  

First note that for $y,z\in X$, $y^{g(y)}\cdot
z^{g(z)}=0$ holds if and only if
there are $\alpha\in S$ and $i\in\omega$ with
$\{y,z\}=\{x_\alpha,x_{\delta_i^\alpha}\}$ such that
$g(x_\alpha)=1$ and $g(x_{\delta_i^\alpha})=0$.
Now we show the claim by in induction on 
$\max\{\alpha<\lambda:x_\alpha\in Y\}$.  The case
$\card{Y}<3$ is trivial.

Assume the claim has been proved for 
$\max\{\alpha<\lambda:x_\alpha\in Y\}<\beta$.  Suppose 
$\max\{\alpha<\lambda:x_\alpha\in Y\}=\beta$ and for no two elements $y,z\in
Y$, $y^{g(y)}\cdot z^{g(z)}=0$.
For $\beta\not\in S$ the argument is easy.  By assumption,
$b:=\prod_{x\in Y\setminus\{x_\beta\}}x^{g(x)}\not=0$.  By construction, 
$x_\beta$ and $b$ are independent.  Thus 
$\prod_{x\in Y}x^{g(x)}\not=0$.

Now suppose $\beta\in S$ and $\prod_{x\in Y}x^{g(x)}=0$.
By construction, $A_\beta\restriction -x_\beta=\{0\}$.
Thus $b:=\prod_{x\in Y\setminus\{x_\beta\}}x^{g(x)}\not\leq -x_\beta$.
Therefore $g(x_\beta)=1$ and $b\leq x_\beta$.  
By construction, there is $m\in\omega$
such that $b\leq\sum_{i<m}x_{\delta_i^\beta}$.   
It follows
from the inductive hypothesis that
$b\cdot\prod_{i<m}-x_{\delta_i^\beta}\not=0$.  This contradicts the choice
of
$m$ and the claim is proved.

For each $a\in C$ let $X_a=\{x_{a,i}:i<n\}$.  Clearly, we may assume that
$C$
has size $\aleph_1$.  Let $\leq$ be a wellordering on $C$ of ordertype
$\omega_1$.  For each $\{a,b\}\in[C]^2$ choose a color $c(\{a,b\})\in n^2$
such that 
$$\forall(i,j)\in n^2(c(\{a,b\})=(i,j)\wedge a\leq b
\Rightarrow x_{a,i}^{f_a(x_{a,i})}\cdot x_{b,j}^{f_b(x_{b,j})}=0).$$

It follows from the claim that $c$ can be defined.
Clearly, for all $\{a,b\}\in[C]^2$, if $c(\{a,b\})=(i,j)$ and $a\leq b$, then
$x_{a,i},x_{b,j}\not\in R$. 
In \cite{baha} Baumgartner and Hajnal established the
following partition
result:
$$\forall
m\in\omega\forall\alpha<\omega_1(\omega_1\rightarrow(\alpha)_m^2).$$

In particular, $\omega_1\rightarrow(\omega+2)^2_{n^2}$ holds.
That is, there are $(i,j)\in n^2$ and a 
subset $C^\prime$ of $C$ of ordertype $\omega+2$
such that for all $\{a,b\}\in[C^\prime]^2$,  $c(\{a,b\})=(i,j)$.
Let $a$ and $b$ be the last two elements of $C^\prime$.
Assume $x_{a,j}=x_\alpha$ for some $\alpha\in S$.  By construction of $A$,
for all $c\in C^\prime\setminus\{a,b\}$, $x_{c,i}=x_{\delta_k^\alpha}$ for
some $k\in\omega$.  By the $\Delta$-system assumption, all the $x_{c,i}$'s 
are
different.  This implies $x_{a,j}=x_{b,j}$, contradicting the
$\Delta$-system assumption.

Now assume that for all $\alpha\in S$, $x_{a,j}\not=x_\alpha$.  In this
case, for all $c\in C^\prime\setminus\{a,b\}$,
$x_{c,i}=x_\alpha$ for some $\alpha\in S$.  Let $d$ and $e$ be the first two
elements of $C^\prime$.  Now for all $c\in C^\prime\setminus\{d,e\}$, 
$x_{c,j}=x_{\delta_k^\alpha}$ for
some $k\in\omega$.  By the $\Delta$-system assumption, all the $x_{c,j}$'s 
are
different.  This implies $x_{d,i}=x_{e,i}$, contradicting the
$\Delta$-system assumption.  This finishes the proof of the lemma.
\end{proof}
 
\begin{proof}[Proof of Theorem \ref{tightmany}] Let $\kappa$ be an
infinite cardinal.
If $\kappa=\aleph_0$, then there are 
$2^\kappa$ pairwise non-isomorphic Boolean algebras of
size $\kappa$ and all of them are projective, 
hence tightly $\sigma$-filtered.
Also, if $\kappa$ is singular, then there are 
$2^\kappa$ pairwise non-isomorphic
projective Boolean algebras by the result of Koppelberg mentioned before. 
Projective Boolean algebras are c.c.c.

If $\kappa$ is regular and uncountable, let $\mathcal P$ 
be a disjoint family of
stationary subsets of $\{\alpha<\kappa:\cf(\alpha)=\aleph_0\}$ of size
$\kappa$.  
Such a family exists by the wellknown results of Ulam and Solovay.  For every
subset $\mathcal T$ of $\mathcal P$ let $A^{\mathcal T}$ be the Boolean
algebra which is constructed in Lemma \ref{tightcoding} from the set
$S:=\bigcup\mathcal T$ and let $(A_\alpha^{\mathcal T})_{\alpha<\kappa}$
be 
its associated tight $\sigma$-filtration.   Then for $\mathcal T,\mathcal
T^\prime\subseteq\mathcal P$ with $\mathcal T\not=\mathcal T^\prime$ the 
Boolean algebras $A^{\mathcal T}$ and $A^{\mathcal T^\prime}$ are 
non-isomorphic. 

For
suppose $h:A^{\mathcal T}\longrightarrow A^{\mathcal T^\prime}$ is an 
isomorphism.
Without loss of generality we may assume that $\mathcal T\setminus\mathcal
T^\prime$ is nonempty.
The set 
$\{\alpha<\kappa:h[A^{\mathcal T}_\alpha]=A^{\mathcal T^\prime}_\alpha\}$
is club in $\lambda$.  
Since $\bigcup(\mathcal T\setminus\mathcal T^\prime)$ 
is stationary, there is
$\alpha\in\bigcup(\mathcal T\setminus\mathcal T^\prime)$ 
such that $h[A_\alpha^{\mathcal T}]=A^{\mathcal T^\prime}_\alpha$.
But $A^{\mathcal T}_\alpha\not\leq_{\rc}A^{\mathcal T}$ and $A^{\mathcal
T^\prime}_\alpha\leq_{\rc}A^{\mathcal T^\prime}$, a contradiction.

By Lemma \ref{ccc}, the Boolean algebras $A^{\mathcal T}$ are
c.c.c.
\end{proof}

The two lemmas above give even more: 
\begin{thm}  Let $\kappa$ be an uncountable and regular cardinal.  
Then there is a family of size
$2^\kappa$ of tightly $\sigma$-filtered c.c.c.~Boolean 
algebras of size $\kappa$ 
such that no member of this family is embeddable into another one as an
rc-subalgebra. 
\end{thm}
        
\begin{proof}
Suppose $\mathcal T$ and $\mathcal T^\prime$ are subsets of $\mathcal P$,
where $P$ is as in
the proof of the theorem above.
Assume there is an embedding 
$e:A^{\mathcal T}\to A^{\mathcal T^\prime}$
such that $e[A^{\mathcal T}]\leq_{\rc}A^{\mathcal T^\prime}$.   

Let $C\subseteq\kappa$ be a club 
such that 
$e[A_\alpha^{\mathcal T}]=A^{\mathcal T^\prime}_\alpha\cap e[A^{\mathcal T}]$ 
and $\lpr^{A^{\mathcal T^\prime}}_{e[A^{\mathcal T}]}[A^{\mathcal
T^\prime}_\alpha]\subseteq A^{\mathcal T^\prime}_\alpha$ hold for
every $\alpha\in C$.

Let $\alpha\in C\cap\bigcup\mathcal T$. Then 
$e[A^{\mathcal T}_\alpha]\not\leq_{\rc}e[A^{\mathcal T}]$ and hence 
$e[A^{\mathcal T}_\alpha]\not\leq_{\rc}A^{\mathcal T^\prime}$.
Since
$A^{\mathcal T^\prime}_\alpha$ is closed under 
$\lpr^{A^{\mathcal T^\prime}}_{e[A^{\mathcal T}]}$,
$e[A_\alpha^{\mathcal T}]\leq_{\rc}A_\alpha^{\mathcal T^\prime}$.

Hence
$A_\alpha^{\mathcal T^\prime}\not\leq_{\rc}A^{\mathcal T^\prime}$.
Therefore $C\cap\bigcup\mathcal T\subseteq C\cap\bigcup\mathcal T^\prime$.
Thus, since $\mathcal P$ consists of stationary sets,
$\mathcal T\subseteq\mathcal T^\prime$.
Now let $I$ be an independent family of subsets of $\mathcal P$ of size
$2^\lambda$.  In particular, the elements of $I$ are pairwise 
$\subseteq$-incomparable.
Thus the family $\{A^{\mathcal T}:\mathcal T\in I\}$ consists of pairwise
non-rc-embeddable tightly $\sigma$-filtered 
c.c.c.~Boolean algebras of size $\lambda$.
\end{proof}

\section{Openly generated Boolean algebras}

\subsection{Almost free Boolean algebras}
The openly generated Boolean algebras we are going to construct will be
almost free.  Unlike in the case of abelian groups, subalgebras of free
Boolean algebras do not have to be free.  Thus we have to use a slightly
generalized definition of almost freeness for Boolean algebras.
We use the definition given in \cite{ekmek}.

\begin{defn} Let $A$ be a Boolean algebra of size $\kappa$.  
A filtration $(A_\alpha)_{\alpha<\kappa}$ of $A$ is a {\em
$\kappa$-filtration} 
if for every $\alpha<\kappa$, $\card{A_\alpha}<\kappa$. 

$A$ is {\em almost free} if it has a $\kappa$-filtration
$(A_\alpha)_{\alpha<\kappa}$ consisting of free Boolean algebras.  
\end{defn}

In the proof of Theorem \ref{tightmany} we coded stationary sets by
tightly $\sigma$-filtered Boolean algebras using the difference
between rc-subalgebras and $\sigma$-subalgebras.  
If there is an almost free abelian group of size $\kappa$, 
we can code certain stationary subsets of $\kappa$ by openly generated
Boolean algebras.  The coding will be more subtle as in the case of
tightly $\sigma$-filtered Boolean algebras.

\begin{defn} Let $A$ and $B$ be Boolean algebras such that $B\leq A$.
Then $B\leq_{\free}A$ if there is a free Boolean algebra $F$ such that $A$
is isomorphic to $B\oplus F$ over $B$.  Here $\oplus$ denotes the
coproduct (free product) in the category of Boolean algebras.
Equivalently, $B\leq_{\free}A$ if there is a set $X\subseteq A$ such that
$X$ is independent over $B$ and $A=\langle B\cup X\rangle$.  

Suppose $A$ is almost free of size $\kappa$ and let
$(A_\alpha)_{\alpha<\kappa}$ be a $\kappa$-filtration of $A$.
Let $$E:=\{\alpha<\kappa:\{\beta\in(\alpha,\kappa):A_\alpha\not\leq_{\free} A_\beta\}\mbox{ is stationary in }\kappa\}.$$
For $X\subseteq\kappa$ let
$$\tilde{X}:=\{Y\subseteq\kappa:\exists
C\subseteq\kappa(C\mbox{ is
club and }Y\cap C=X\cap C)\}.$$ 
Let $\Gamma(A):=\tilde{E}$.
For $X,Y\subseteq\kappa$ let
$\tilde{X}\leq\tilde{Y}$ if $X\setminus Y$ is non-stationary.  
\end{defn}

It is routine matter to check that $\Gamma(A)$ does not depend on the
choice of $(A_\alpha)_{\alpha<\kappa}$.  Therefore
$\Gamma(A)=\Gamma(B)$ if $A$ and $B$ are isomorphic.
Thus $\Gamma(A)$ is an invariant of $A$.  
Also, it is not difficult to see that $A$ is free if and only if
$\Gamma(A)=\tilde{\emptyset}$.  See \cite{ekmek} for more information on
the $\Gamma$-invariant.

\subsection{Almost free families of countable sets and the strong
construction principle $\operatorname{(CP+)}$}

Shelah has given an exact translation of the algebraic question whether
there is an
almost free, non-free abelian group of size $\kappa$ into a set-theoretic
one.  We use the representation given in \cite{ekmek}, where all the
missing proofs can be found.     

Let $\mathcal S=\{s_i:i\in\mathcal S\}$ be a family of countably
infinite sets.   
A transversal for $\mathcal S$ is a one-one function
$T:I\to\bigcup\mathcal S$ such that for all $i\in I$, $T(i)\in s_i$. 
We say that $\mathcal S$ is {\em free} if it has a transversal.  
$\mathcal S$ is {\em almost free} if for each $J\subseteq I$ with
$\card{J}<\card{I}$, $\{s_i:i\in J\}$ has a transversal.  
Let $\NPT(\lambda,\aleph_0)$ denote the statement
``there is an almost free, non-free family of size $\lambda$ of countable
sets''.  
(We use the notation of \cite{shevai}.)

Shelah proved the following:
\begin{thm}\label{almostfreegroups}  For every uncountable cardinal
$\lambda$ there is an almost
free, non-free abelian group of size $\lambda$ if and only if $\NPT(\lambda,\aleph_0)$ holds.
\end{thm}

We mimic the proof of one direction of this theorem
and construct a Boolean algebra $A(\mathcal S)$ for any given family
$\mathcal S$ of
countable sets.  For this we use the strong construction
principle $\operatorname{(CP+)}$ for Boolean algebras.  If the family $\mathcal S$ is
sufficiently good, $A(\mathcal S)$ will be openly generated.

$\operatorname{(CP+)}$ (for arbitrary varieties) was introduced 
by Eklof and Mekler, who  proved

\begin{thm}\label{varietywithCP+} $\NPT(\lambda,\aleph_0)$ implies the existence of
almost free, non-free objects of size $\lambda$ in every variety
$\mathcal V$ which satisfies $\operatorname{(CP+)}$.
\end{thm} 

Shelah showed 
\begin{lemma}\label{CP+}(See \cite{fuhabil}) 
The variety of Boolean algebras satisfies the
following {\em strong construction principle} $\operatorname{(CP+)}$:

For each $n\in\omega\setminus 1$ there are countably generated free
Boolean algebras 
$H\leq K\leq L$ and a partition of $\omega$ into $n$ infinite blocks 
$s^1,\dots,s^n$ such that 
\begin{itemize}\item[(i)] $H$ is freely generated by $\{h_m:m\in\omega\}$
and for each $J\subseteq\omega$, if for some $k\in\{1,\dots,n\}$, $J\cap
s^k$ is finite, then $\langle\{h_m:m\in J\}\rangle\leq_{\free} L$ and
\item[(ii)] $L=K\oplus\Fr(\omega)$ and $H\not\leq_{\free} L$.  
\end{itemize}
\end{lemma}  

For the convenience of the reader we include a proof of this lemma.  
The proof uses
\begin{lemma}(Sirota's Lemma, see \cite[Theorem 1.4.10]{heinsha})\label{sirota}
Let $A$ and $B$ be Boolean algebras such that $A\leq_{\rc} B$ and $B$ is countably generated over $A$.  If for all $b_1,\dots,b_n\in B$ there is $u\in B$ such that $u$
is independent over $A(b_1,\dots,b_n)$, then $A\leq_{\free}B$.  
\end{lemma}

\begin{proof}[Proof of Lemma \ref{CP+}] Let $n\in\omega\setminus 1$ and let
$(x_{k,l})_{k\in\{1,\dots,n\},l\in\omega}$ be a family of pairwise distinct sets.
Let $H$ be the free Boolean algebra $\Fr(X)$ over the set $X:=\{x_{k,l}:k\in\{1,\dots,n\}\wedge l\in\omega\}$.
Let $(h_m)_{m\in\omega}$ be a 1-1-enumeration of $X$.  
For $k\in\{1,\dots,n\}$ let $s^k:=\{m\in\omega:\exists l\in\omega(h_m=x_{k,l})\}$.
Let $I$ be the ideal of $H$ generated by $\{\prod_{l<n}x_{l,h}:h\in\omega\}$.  
Let $K$ be a Boolean algebra of the form $H(x)$ where $H\restriction x=I$ and $H\restriction-x=\{0\}$.   
Finally, let $L:=K\oplus\Fr(\omega)$.  

We show that $H$, $K$, and $L$ are as required in the definition of $\operatorname{(CP+)}$.
Clearly, the three Boolean algebras are countable and atomless.  
Therefore, they are free.
By the choice of $K$ and $x$, $H\not\leq_{\rc}K$ and thus $H\not\leq_{\rc}L$ .  
In particular, $H\not\leq_{\free} L$.  

Now suppose that $Y\subseteq X$ is such that for some $k_0\in\{1,\dots,n\}$ 
the set $Y\cap\{x_{k_0,l}:l\in\omega\}$ is finite.  
Let $H^\prime$ be the subalgebra of $H$ generated by $Y$.    
We have to show $H^\prime\leq_{\free}L$.  
By Lemma \ref{sirota} it is enough to show $H^\prime\leq_{\rc}L$.  
Since $K\leq_{\free}L$ and thus $K\leq_{\rc}L$, it is in fact sufficient 
to show $H^\prime\leq_{\rc}K$.

For every $m\in\omega$ let $H_m$ be the subalgebra of $H$ generated by 
$$\{x_{k,l}:k\in\{1,\dots,n\}\wedge l\in\omega\wedge(k=k_0\Rightarrow l\leq m)\}.$$
Let $a\in K=H(x)$.  
Then for some $m\in\omega$, $H^\prime\subseteq H_m$ and $a\in H_m(x)$.
Clearly, $H^\prime\leq_{\rc}H_m$.
Since $I\cap H_m=\{0\}$, $x$ is independent over $H_m$.
Therefore $H_m\leq_{\rc}H_m(x)$.  
It follows that $H^\prime\leq_{\rc}H_m(x)$.  
This implies that $H^\prime\restriction a$ is principal.  
This shows $H^\prime\leq_{\rc}H(x)$, finishing the proof of the lemma.  
\end{proof}

The construction of $A(\mathcal S)$ from a
family $\mathcal S$ of countable sets is the following.

\begin{defn}\label{BAfromS}  Let $\mathcal S=\{s_i:i\in I\}$ be a family
of
countable sets.  
Suppose for all $i\in I$, $s_i$ is the disjoint union of the infinite
sets $s_i^1,\dots,s_i^n$. Let $H$, $K$, $L$, $\{h_m:m\in\omega\}$, and
$s^1,\dots,s^n$ be as in Lemma \ref{CP+}.
For each $i\in I$ fix an enumeration $(x_{i,m})_{m\in\omega}$ of $s_i$
such that  for each $k\in\{1,\dots,n\}$, $s_i^k=\{x_{i,m}:m\in s^k\}$.  

For each $i\in I$ choose a copy $L_i$ of $L$.  Assume that all the copies
are disjoint. 
Let $K_i$, $H_i$, and $\{h_{i,m}:m\in\omega\}$ denote the corresponding
copies of $K$, $H$, and $\{h_m:m\in\omega\}$ in $L_i$.  

Let $G:=\bigoplus_{i\in I}L_i$ and let
$\Theta$ be the smallest congruence on $G$ identifying $h_{i,m}$
and
$h_{j,l}$ for all $i,m,j,l$ with $x_{i,m}=x_{j,l}$.  Let
$A(\mathcal S):=G/\Theta$.  
\end{defn}

\subsection{$\lambda$-systems and the strong reshuffling property}
We can only show that the algebra $A(\mathcal S)$ constructed in
Definition \ref{BAfromS} is openly generated if the family $\mathcal S$
has some special properties.  We need $\mathcal S$ to be {\em based on 
a $\lambda$-system $\Lambda$ of height $n$} for some $n>0$ such that
$(\mathcal S,\Lambda)$ has the {\em strong reshuffling property}.  

Shelah used $\lambda$-systems and the (weak) reshuffling property for
proving Theorem \ref{almostfreegroups}.   The proof of Theorem
\ref{varietywithCP+} uses the strong reshuffling property. 
We use $\lambda$-systems to construct many openly
generated Boolean algebras in the same way as 
they are used in the
proof of
Theorem \ref{varietywithCP+}.
Let us start with the definition of a $\lambda$-system.

\begin{defn} The set $\lambda^{<\omega}$ ordered by set inclusion is a
tree.  Let $(\beta)$ denote the sequence of length one with value $\beta$ 
and let $\frow$ denote concatenation of sequences.  
If $S$ is a subtree of $\lambda^{<\omega}$, an element $\eta$ of $S$ is
called a {\em final node} of $S$ if in $S$ there is no proper extension
of $\eta$.  Let $S_f$ denote the set of final nodes of $S$.

1. A {\em $\lambda$-set} $S$ is a subtree of $\lambda^{<\omega}$
together
with a cardinal $\lambda_\eta$ for every $\eta\in S$ such that
$\lambda_\emptyset=\lambda$ and
\begin{itemize}\item[a)] for all $\eta\in S$, $\eta\in S_f$ if and only if
$\lambda_\eta=\aleph_0$ and 
\item[b)] if $\eta\in S\setminus S_f$, then $\eta\frow(\beta)\in S$
implies $\beta\in\lambda_\eta$ and $\lambda_{\eta\frow
(\beta)}<\lambda_\eta$ and
$E_\eta:=\{\beta<\lambda_\eta:\eta\frow(\beta)\in S\}$ is stationary in
$\lambda_\eta$.
\end{itemize}

2. A {\em $\lambda$-system} is a $\lambda$-set together with a set
$B_\eta$
for each $\eta\in S$ such that $B_\emptyset=\emptyset$ and for all
$\eta\in S\setminus S_f$
\begin{itemize}
\item[a)] for all $\beta\in E_\eta$,
$\lambda_{\eta\frow(\beta)}\leq\card{B_\eta}<\lambda_\eta$ and
\item[b)] $(B_{\eta\frow(\beta)})_{\beta\in E_\eta}$ is increasing and
continuous, that is, if $\sigma\in E_\eta$ is a limit point of $E_\eta$,
then $B_{\eta\frow(\sigma)}=\bigcup\{B_{\eta\frow(\beta)}:\beta\in
E_\eta\cap\sigma\}.$
\end{itemize}

For any $\lambda$-system $\Lambda=(S,\lambda_\eta,B_\eta)_{\eta\in S}$ and
any $\eta\in S$ let $\overline{B}_\eta:=\bigcup\{B_{\eta\restriction
m}:m\leq\dom(\eta)\}$.  A family $\mathcal S$ of countable sets is {\em
based on} $\Lambda$ if $\mathcal S$ is indexed by $S_f$ and for every
$\eta\in S_f$, $s_\eta\subseteq\overline{B}_\eta$.  

A subtree $S$ of $\lambda^{<\omega}$ has {\em height} $n$ if all final
nodes of $S$ have domain $n$.  A $\lambda$-set or $\lambda$-system has
height $n$ if its associated tree $S$ has height $n$.  It is not difficult
to see that every $\lambda$-set has a sub-$\lambda$-set which has a
height.  
\end{defn}
  
If $\mathcal S$ is a family of countable sets based on a $\lambda$-system,
then $\mathcal S$ has cardinality $\lambda$ and is not free.
Families of countable sets based on $\lambda$-systems can be constructed
from almost free, non-free abelian groups (see \cite{ekmek}).  
By Theorem \ref{almostfreegroups}, this implies

\begin{lemma}\label{lambdasystem} $\NPT(\lambda,\aleph_0)$ implies that 
there is a family of countable sets based on a $\lambda$-system.  
\end{lemma}

In \cite{shevai} Shelah and V\"ais\"anen provided the tools to show that 
for $\lambda>\aleph_1$, the $\lambda$-system in Lemma \ref{lambdasystem} 
may be chosen such that its height is at
least $2$, i.e., such that the underlying $\lambda$-set is really more than just a
stationary subset of $\lambda$.

\begin{lemma}\label{increasingheight}
If $\lambda>\aleph_1$ and $\NPT(\lambda,\aleph_0)$ holds, 
then there is a family of countable sets based on a
$\lambda$-system of height $n$ for some $n>1$.
\end{lemma}

\begin{proof}
We show how to extract the proof of Lemma \ref{increasingheight} 
from \cite{shevai}\footnote{The second author thanks Pauli V\"ais\"anen
for explaining how Lemma \ref{increasingheight} follows from 
the results presented in \cite{shevai}.}.  
In \cite{shevai} a special kind of families of countable sets based on
$\lambda$-systems is used, so-called $\NPT(\lambda,\aleph_0)$-skeletons.  
Moreover, so-called $\NRT(\lambda,\aleph_0)$-skeletons are used for building more
complicated $\NPT$-skeletons.  

Let $\lambda>\aleph_1$ and suppose $\NPT(\lambda,\aleph_0)$ holds.  
We show that there is an $\NPT(\lambda,\aleph_0)$-skeleton of height $>1$.  

\relax From an $\NRT(\lambda,\aleph_0)$-skeleton one can construct an 
$\NPT(\lambda,\aleph_0)$-skeleton of the same height \cite[Corollary 4.14]{shevai}.  
It is therefore sufficient to show the existence of an 
$\NRT(\lambda,\aleph_0)$-skeleton of height $>1$.  
By $\NPT(\lambda,\aleph_0)$, there is an $\NRT(\lambda,\aleph_0)$-skeleton
\cite[Lemma 4.13 and Corollary 4.14]{shevai}.  
If this $\NRT(\lambda,\aleph_0)$-skeleton is of height $>1$, we are done.  
So assume it is of height $1$.  
Since the $\NRT(\lambda,\aleph_0)$-skeleton is of height
$1$, its type is $(\lambda)$.

As mentioned in the introduction, there is an almost free, non-free abelian group of
size $\aleph_1$.  
By Theorem \ref{almostfreegroups}, this implies $\NPT(\aleph_1,\aleph_0)$.  
As before, it follows that there is an $\NRT(\aleph_1,\aleph_0)$-skeleton. 
An $\NRT(\aleph_1,\aleph_0)$-skeleton is of height $1$. 
Since the $\NRT(\lambda,\aleph_0)$-skeleton is of type $(\lambda)$ and since 
$\lambda>\aleph_1$, the $\NRT(\lambda,\aleph_0)$-skeleton and the 
$\NRT(\aleph_1,\aleph_0)$-skeleton are compatible.  
Therefore, the two skeletons can be combined to an 
$\NRT(\lambda,\aleph_0)$-skeleton of height $2$ (and of type $(\lambda,\aleph_1)$).  \end{proof}

To get an almost free Boolean algebra from $\mathcal S$, we need $\mathcal
S$ to be based on a $\lambda$-system such that $(\mathcal S,\Lambda)$
has the strong the reshuffling property.  

\begin{defn} Let $\Lambda=(S,\lambda_\eta,B_\eta)_{\eta\in S}$ be a
$\lambda$-system of height $n$ and $\mathcal S=(s_\eta)_{\eta\in S_f}$ a
family of
countable sets based on $\Lambda$.  
Let $<_{\lex}$ denote the (strict)
lexicographic order on $S$.

$(\mathcal S,\Lambda)$ has the {\em strong reshuffling property}
if the following statements hold:

\begin{enumerate}
\item If $I$ is a subset of $S_f$ of cardinality $<\lambda$ and $\eta_0\in I$,
then there is a wellordering $<_I$ of $I$ such that $\eta_0$ is the 
$<_I$-first element of $I$ and for all $\eta\in I$ there is $k\in\{1,\dots,n\}$
such that $s_\eta^k\cap\bigcup_{\nu<_I\eta}s_\nu$ is finite.


\item 
For all $\mu\in S\setminus S_f$ and all $\alpha<\lambda_\mu$, 
if $I$ is a subset of 
$\{\eta\in S_f:\mu\subseteq\eta\}$ 
of size $<\lambda_\mu$, then there is a wellordering 
$<_{\bar I}$ of $\bar I:=\{\eta\in S_f:\eta<_{\lex}\mu\vee\eta\in I\}$
such that for all $\eta\in\bar I$
\begin{itemize}\item[a)]
there is $k\in\{1,\dots,n\}$ such that
$s_\eta^k\cap\bigcup_{\nu<_{\bar I}\eta}s_\nu$ is finite and
\item[b)] 
if $\mu\subsetneqq\eta$ and 
$\nu\in\bar I$ is such that
$$\nu<_{\lex}\mu\vee(\mu\subsetneqq\nu\wedge\nu(\dom(\mu))
\leq\alpha<\eta(\dom(\mu))),$$
then $\nu<_{\bar I}\eta$.
\end{itemize}

\end{enumerate}
\end{defn}

The definition of the strong reshuffling property given in
\cite{ekmek} is slightly weaker than ours.   
However, even with our version of the strong reshuffling
property, the proof of Theorem 3A.7  in \cite{ekmek} shows

\begin{lemma}\label{reshuffling}
If there is a family of countable sets based on a $\lambda$-system of
height $n$, then 
there is a family $\mathcal S$ of countable sets based on a
$\lambda$-system $\Lambda$ of height $n$ such that
$(\mathcal S,\Lambda)$ has the strong reshuffling property.  
\end{lemma}

Combining this with Lemma
\ref{increasingheight} we get

\begin{corollary}\label{goodlambdasystem} If $\lambda>\aleph_1$ and 
$\NPT(\lambda,\aleph_0)$ holds, then there is a family $\mathcal
S$ of countable sets based on a $\lambda$-system $\Lambda$ of height $>1$ 
such that $(\mathcal S,\Lambda)$ has the strong reshuffling
property.
\end{corollary}

\subsection{Many openly generated Boolean algebras}

In \cite[Theorem 3A.13]{ekmek} it is proved that the algebra $A(\mathcal S)$ 
is almost free but not free if the family $\mathcal S$ of countable sets
is based on a $\lambda$-system $\Lambda=(S,\lambda_\eta,B_\eta)_{\eta\in S}$ 
such that $(\mathcal S,\Lambda)$ has the
strong reshuffling property.  
Looking at the proof more closely, it turns
out that we actually get that $\Gamma(A(\mathcal S))$ is
$\tilde{E}_\emptyset$.
In the case of Boolean algebras we find that $A(\mathcal S)$ is
openly generated if $\Lambda$ has height $>1$.
The proof of open generatedness is where we need our stronger version
of the strong reshuffling property.  

\begin{lemma}\label{codingopen}
Let $\mathcal S$ be a family of countable sets based on a
$\lambda$-system 
$\Lambda=(S,\lambda_\eta,B_\eta)_{\eta\in S}$ of height $n>1$ such that
$(\mathcal S,\Lambda)$ has the strong
reshuffling property. 
Then  $\Gamma(A(\mathcal S))=\tilde{E}_\emptyset$ 
and the Boolean algebra $A(\mathcal S)$ is almost free and openly generated.
\end{lemma}

\begin{proof} We only have to show that $A(\mathcal S)$ is openly generated.  
For this we repeat the proof that $A(\mathcal S)$ is almost free.  

Consider the $\lambda$-filtration $(A_\alpha)_{\alpha<\lambda}$
of $A(\mathcal S)$ where
$$A_\alpha:=\langle\bigcup\{L_\eta:\eta\in
S_f\wedge\eta(0)<\alpha\rangle/\Theta$$
for all $\alpha<\lambda$.

{\bf Claim 1.} $A_{\alpha+1}\leq_{\free}A_\beta$ for all
$\beta<\lambda$ and all $\alpha\in[-1,\beta)$.  

Note that $A_\alpha=A_{\alpha+1}$ for $\alpha\not\in E_\emptyset$. 
Thus, Claim 1 implies that
\begin{itemize}\item[(i)] all $A_\alpha$, $\alpha<\lambda$, are free and
\item[(ii)] $A_\alpha\leq_{\rc} A(\mathcal S)$ for all $\alpha\not\in
E_\emptyset$.
\end{itemize}

For the proof of Claim 1 let $\alpha<\beta<\lambda$. 
Let $I:=\{\eta\in S_f:\eta(0)<\beta\}$. 
By the strong reshuffling property, there is 
a wellordering $<_I$ of $I$ as guaranteed by the requirement 
(2) in the definition of the strong reshuffling property.  

Now let $\eta$ be the $<_I$-first element of $I$ such that
$\eta(0)>\alpha$.   For some $k\in\{1,\dots,n\}$,
$s_\eta^k\cap\bigcup_{\nu<_I\eta}s_\nu$ is finite.  
Note that $L_\eta/\Theta\cap A_{\alpha+1}=\langle
\{h_{\eta,m}:x_{\eta,m}\in s_\eta\cap\bigcup_{\nu<_I\eta}
s_\nu\}\rangle/\Theta$. 
Choose $G_\eta\leq L_\eta$ free such that
$$L_\eta=\langle\{h_{\eta,m}:x_{\eta,m}\in s_\eta\cap\bigcup_{\nu<_I\eta}
s_\nu\}\rangle\oplus G_\eta.$$

Then $\langle L_\eta/\Theta\cup A_{\alpha+1}\rangle=A_{\alpha+1}\oplus
G_\eta/\Theta$ and $\Theta$ does not affect $G_\eta$, i.e., the
epimorphism $\pi:G_\eta\to G_\eta/\Theta$ is an isomorphism.  
Thus, we may think of $G_\eta$ as a subalgebra of $A_\beta$. 
By recursion on $<_I$ for all $\tau\in I$ with $\eta<_I\tau$
we can choose a free Boolean
algebra $G_\tau\leq A_\beta$ such that
$A_\beta=A_{\alpha+1}\oplus\bigoplus_{\eta\leq_I\tau}G_\tau$.
This shows Claim 1.

By Lemma \ref{provingopen} 
it remains to show $A_\alpha\leq_{\rc}A(\mathcal S)$ for all
$\alpha\in E_\emptyset$.
Fix $\alpha\in E_\emptyset$. 
Since $A_{\alpha+1}\leq_{\rc}A(\mathcal S)$ and by the 
transitivity of $\leq_{\rc}$, 
$A_\alpha\leq_{\rc}A(\mathcal S)$ if and only if $A_\alpha\leq_{\rc}A_{\alpha+1}$.
Thus, it remains to show $A_\alpha\leq_{\rc}A_{\alpha+1}$.

Since $\Lambda$ is of height $>1$, $(\alpha)\not\in S_f$.  
Consider the filtration $(A_{\alpha,\beta})_{\beta<\lambda_{(\alpha)}}$ of 
$A_{\alpha+1}$ where
$$A_{\alpha,\beta}:=\langle\bigcup\{L_\eta:\eta\in
S_f\wedge\eta(0)<\alpha\vee(\eta(0)=\alpha\wedge\eta(1)<\beta)\rangle/\Theta.$$

{\bf Claim 2.}  $A_\alpha=A_{\alpha,0}\leq_{\free}A_{\alpha,\beta}$ for all
$\beta<\lambda_{(\alpha)}$.  

The proof of Claim 2 is practically the same as the proof of Claim 1 and uses  
requirement (2) in the definition of the strong reshuffling property 
with $\mu:=(\alpha)$ and $I:=\{\eta\in S_f:\eta(0)=\alpha\wedge\eta(1)<\beta\}$.

Now let $a\in A_{\alpha+1}$.  There is $\beta<\lambda_\mu$ such that $a\in A_{\alpha,\beta}$. 
By Claim 2, $A_\alpha\leq_{\free}A_{\alpha,\beta}$.
In particular, $A_\alpha\leq_{\rc}A_{\alpha,\beta}$. 
Therefore $A_\alpha\restriction a$ has a maximal element.  
This shows $A_\alpha\leq_{\rc}A_{\alpha+1}$ and finishes the proof of the lemma.
\end{proof}

Lemma \ref{codingopen} says that in certain cases we can code stationary
sets by openly generated Boolean algebras.  
Using this coding, we can show

\begin{thm}  Let $\lambda>\aleph_1$ be such that there is an almost free,
non-free abelian group of size $\lambda$.  Then there are $2^\lambda$
pairwise non-isomorphic openly generated Boolean algebras of size
$\lambda$.  
\end{thm}

\begin{proof}  If there is an almost free, non-free abelian group of size
$\lambda$, then, by Corollary \ref{goodlambdasystem}, there is a family
$\mathcal S$ of countable sets based on a $\lambda$-system
$\Lambda=(S,\lambda_\eta,B_\eta)_{\eta\in S}$ of height $n>0$ such that
$(\mathcal S,\Lambda)$ has the strong reshuffling property.
By Shelah's compactness theorem, $\lambda$ is regular.  

By the theorems of Solovay and Ulam, we can split the stationary set
$E_\emptyset$ associated with $S$ into a disjoint family $\mathcal
P$ 
of size $\lambda$ of stationary sets.
For each 
$\mathcal T\subseteq\mathcal P$ let
$E^{\mathcal T}:=\bigcup\mathcal T$,
$S^{\mathcal T}:=\{\eta\in S:\eta(0)\in E^{\mathcal T}\}$,
$\Lambda^{\mathcal T}:=
(S^{\mathcal T},\lambda_\eta,B_\eta)_{\eta\in S^{\mathcal T}}$,
$\mathcal S^{\mathcal T}:=(s_\eta)_{\eta\in S^{\mathcal T}_f}$, and 
$A^{\mathcal T}:=A(\mathcal S^{\mathcal T})$.  

It follows immediately from the definitions that for all
${\mathcal T}\subseteq\mathcal P$ with $\mathcal T\not=\emptyset$,
$\Lambda^{\mathcal T}$ is a
$\lambda$-system of height $n$ and 
$(\mathcal S^{\mathcal T},\Lambda^{\mathcal T})$
has the strong reshuffling property.  
Thus, by Lemma \ref{codingopen}, $A^{\mathcal T}$ is openly
generated and $\Gamma(A^{\mathcal
T})=\tilde{E^{\mathcal T}}$ for each nonempty 
${\mathcal T}\subseteq\mathcal P$.  Clearly, 
$\tilde{E^{\mathcal T}}\not=\tilde{E^{\mathcal T^\prime}}$ for 
$\mathcal T\not=\mathcal T^\prime$.  
It follows that $(A^{\mathcal T})_{{\mathcal
T}\subseteq\mathcal P,{\mathcal
T}\not=\emptyset}$ is a family of $2^\lambda$ pairwise non-isomorphic
openly generated Boolean algebras of size $\lambda$.
\end{proof}

Combining this theorem with the results of Shelah and Magidor from
\cite{mashe} on the class of cardinals $\kappa$ for which there are almost
free, non-free abelian groups of size $\kappa$, we get

\begin{corollary} a) There is a class $\mathcal C$ of regular cardinals
which is closed under the operations $\kappa\mapsto\kappa^+$ and
$(\kappa,\lambda)\mapsto\lambda^{+(\kappa+1)}$ and contains $\aleph_2$ 
such that for each $\kappa\in\mathcal C$ there are $2^\lambda$ isomorphism
types of openly generated Boolean algebras of size $\lambda$.   
In particular, for every $n>1$ there are $2^{\aleph_n}$ pairwise
non-isomorphic Boolean algebras of size $\aleph_n$.  

b) Under $V=L$, for every cardinal $\kappa>\aleph_1$ there are
$2^\kappa$
pairwise non-isomorphic openly generated Boolean algebras of size
$\kappa$. 
\end{corollary}

Note that we have to exclude $\aleph_1$ in this corollary since every 
openly generated Boolean algebra of size $\aleph_1$ is projective and
thus, by Koppelberg's result on the number of projective Boolean
algebras, there are only $2^{\aleph_0}$ isomorphism types of 
openly generated Boolean algebras of size $\aleph_1$.

Recall that
open generatedness follows from projectivity.
By Koppelberg's results on the number of projective Boolean algebras of
singular cardinality, for each singular cardinal $\mu$ there are $2^\mu$
isomorphism types of openly generated Boolean algebras of size $\mu$.  

The main open question is

\begin{question} Is it true that the number of isomorphism types of openly
generated Boolean algebras of size $\kappa$ is $2^\kappa$ for every
infinite cardinal $\kappa\not=\aleph_1$?
\end{question}

\end{document}